\documentclass[11pt,onecolumn,twoside]{IEEEtran}

\usepackage{amsmath,amssymb,euscript,yfonts,psfrag,latexsym,dsfont,graphicx,bbm,color,amstext,wasysym,pdfsync,thumbpdf}

\newtheorem{thm}{Theorem}

\newtheorem{prop}[thm]{Proposition}

\newcommand{\mR}{{\mathbb R}}

\newcommand{\cL}{{\mathcal L}}
\newcommand{\cT}{{\mathcal T}}

\newcommand{\tH}{{\text H}}

\newcommand{\tL}{{\text L}}
\newcommand{\tLog}{{\text {Log}}}
\newcommand{\tKL}{{\text {KL}}}

\newcommand{\trace}{{\operatorname{trace}}}
\newcommand{\rank}{{\operatorname{rank}}}

\newcommand{\g}{{\operatorname{g}}}

\begin{document}
\title{Geometric methods for estimation\\ of structured covariances}
\author{Lipeng Ning, Xianhua Jiang and Tryphon Georgiou}
\date{\today}
\maketitle

\begin{abstract}
We consider problems of estimation of structured covariance matrices, and in particular of matrices with a Toeplitz structure. We follow a geometric viewpoint that is based on some suitable notion of distance.
To this end, we overview and compare several alternatives metrics and divergence measures.
We advocate a specific one which represents the Wasserstein distance between the corresponding Gaussians distributions and show that it coincides with the so-called Bures/Hellinger distance between covariance matrices as well. Most importantly, besides the physically appealing interpretation, computation of the metric requires solving a linear matrix inequality (LMI). As a consequence, computations scale nicely for problems involving large covariance matrices, and linear prior constraints on the covariance structure are easy to handle.
We compare this {\em transportation/Bures/Hellinger} metric with the maximum likelihood and the Burg methods as to their performance with regard to estimation of power spectra with spectral lines
on a representative case study from the literature.\footnote{
Dept. of Electrical \& Comp.\ Eng., University of Minnesota, Minneapolis, MN 55455.
{\{ningx015, jiang082, tryphon\}@umn.edu}.\\
\hspace*{12pt}Supported in part by the NSF and AFOSR.}
\end{abstract}

\section{Introduction}

Consider a zero-mean, real-valued, discrete-time stationary random process $\{ x(t),~t\in \mathbb{Z}\}$. Let
\[
r(t):=E\left( x(k)x(k+t)\right),
 \]
with $k,t\in \mathbb{Z}$, denote the autocorrelation function, and
\[
T:=\left[
     \begin{array}{cccc}
       r_0 & r_1 & \cdots & r_{n-1} \\
       r_{-1} & r_0 & \cdots & r_{n-2} \\
       \vdots & \vdots & \ddots & \vdots \\
       r_{-(n-1)} & r_{-(n-2)} & \ldots & r_0 \\
     \end{array}
   \right]
\]
the covariance of the finite (observation) vector
\newcommand{\bx}{{\mathbf x}}
\newcommand{\bX}{{\mathbf X}}
 \[
\bx= \left[\begin{matrix}x(0),&x(1),&\ldots&x(n-1)\end{matrix}\right]^\prime,
 \]
 i.e., $T=E(\bx\bx^\prime)$.
The covariance has a Toeplitz structure inherited by the time-invariance (stationarity) of the process.
Throughout, the size of such an observation vector and of corresponding finite Toeplitz matrices will always be $n$ and $n\times n$, respectively.

The power spectrum of the process is uniquely determined by the (infinite) autocorrelation function. This is due to the fact that the trigonometric moment problem is {\em determined} \cite{grenander2001toeplitz}. Then, starting with Burg's early contributions \cite{Burg_thesis,haykin1979nonlinear}, modern nonlinear spectral analysis techniques largely rely on {\em admissible} estimates of the partial autocorrelation sequence $\{r_0,\;r_1,\ldots,r_{n-1}\}$ (equivalently, of the Toeplitz covariance $T$) from which information is sought about corresponding power spectra. Admissibility of the partial autocorrelation sequence amounts to the requirement that $T$ is a positive semi-definite matrix, in which case a positive semi-definite extension to an infinite matrix is also possible.

Part of the challenge, which was already addressed by Burg, is due to the fact that the sample covariance
\begin{align}\label{Tsample}
\hat{T} :=\frac{1}{m}\sum_{k=1}^m \bx_k\bx_k^\prime,
\end{align}
where $\bx_k$ ($k\in\{1,2,\ldots,m\}$) are independent observation vectors,
{\em may not be Toeplitz} due to statistical errors. On the other hand, estimates of the individual entries
$\{r_0,\;r_1,\ldots,r_{n-1}\}$ via averaging
over all available samples to within a given time-distance from one another, may not lead to
a positive matrix. Either way, the linear structure or the positivity is compromised.

An early popular algorithm by Burg aimed at ensuring positivity via a clever estimation of the so-called partial reflection coefficients instead of the autocorrelation coefficients (see e.g., \cite{haykin1979nonlinear}). Several alternative tricks were devised followed by a maximum likelihood approach in \cite{BLW_estimation}.
However, the issue was never put to rest because all these face challenges of their own that lead to poor resolution, bias, ``line-spliting'' (where sinusoidal components in the spectrum generate ghost peaks), and computational difficulties (as in the case of \cite{BLW_estimation}).
The source of the problem is largely the error in $T$ which adversely affects our subsequent estimate of the underlying power spectrum (obtained using e.g., a Maximum Entropy method, the Capon envelope, etc.). Herein, we do not analyze the problem of going from the Toeplitz covariance to a power spectral estimate. Instead we focus only on the problem of estimating the Toeplitz covariance from finite observations.

The Toeplitz covariance matrix is sought as the one closest to $\hat{T}$ in a suitable geometry. Notions of distance from information theory, quantum mechanics, and statistics lead to complementary viewpoints and so does maximum likelihood estimation of the autocorrelation coefficients which also provides us with a notion of distance.
In Section \ref{sec:geometric} we outline the geometric viewpoint together with various possibilities for distance measures. In Section \ref{sec:approximation} we discuss the respective optimization problems and in Section \ref{sec:examples} we compare the three most promising alternatives on a specific example from the literature.

\section{Geometric Viewpoint}\label{sec:geometric}
Given a sample covariance matrix $\hat T$, we consider the problem to minimize
\begin{align}\label{approximation}
\min_{T\in\cT} d(T,\hat{T} ),
\end{align}
over the class of admissible matrices
\[
    \cT:=\{T: T\geq 0,~T~\text{being Toeplitz} \}.
\]
In this, $d$ represents a suitable notion of distance.
Various such distance measures are motivated below based on statistics, information theory, quantum mechanics, and optimal transportation. Occasionally, when the distance measure is not symmetric, we use the notation $d(T\| \hat{T} )$ instead. Such non-symmetric measures are often referred to as divergences in the literature.

\subsection{Likelihood divergence}
We begin by discussing maximum likelihood estimation \cite{BLW_estimation}.
Assuming that the process $\left\{ x(t), t\in{\mathbb Z} \right\}$ is Gaussian, the joint density function for independent observation vectors $\bx_k$ ($k\in\{1,2,\ldots,m\}$) is
\begin{align*}
p(\bX;T)=(2\pi)^{-\frac{mn}{2}}|T|^{-\frac{m}{2}}\exp\left(-\frac12\sum_{k=1}^{m}\bx_k' T^{-1} \bx_k\right),
\end{align*}
with $\bX:=\left[\bx_1,\ldots,\bx_m\right]$. Then, $\hat T= \frac{1}{m}\bX\bX^\prime$ and
the log-likelihood function becomes
\begin{align}
\cL(\hat{T} ,T)&=\log p(\bX; T )\nonumber\\
&=-\frac{m}{2}\left( n\log(2\pi)+\log|T|+\trace (\hat{T} T^{-1})\right).\label{loglikelihood}
\end{align}
Thus, it is natural to seek a Toeplitz covariance matrix $T$ for which $\cL(\hat{T} ,T)$ is maximal.
Note that if the $\bx_k$'s are independent Gaussian random variables,  $m\hat{T}$ follows a Wishart distribution. Then \eqref{loglikelihood} is the log-likelihood function of this distribution.

Alternatively, one may consider the likelihood divergence
\begin{align}\label{dist_likelihood}
d_{\tL}(T ||\hat{T} )&:=\frac{1}{m}(\log p(\bx; \hat{T} )-\log p(\bx; T) )\nonumber\\
&=\frac{1}{2}(-\log|\hat{T} |+\log|T|+\trace (\hat{T} T^{-1})-n) \nonumber
\end{align}
as a relevant notion of distance since, evidently,
\begin{align*}
&d_{\tL}(T ||\hat{T} )\geq0,\\
&d_{\tL}(T ||\hat{T} )=0\Leftrightarrow T=\hat{T}.
\end{align*}
It relates to the Kullback-Leibler divergence between corresponding pdf's which is discussed next.
However, it does not define a metric because it lacks symmetry and may also fail to satisfy the triangular inequality.

\subsection{Kullback-Leibler divergence}
For random variables on $\mR^n$ with probability density functions $p$ and $\hat p$, the Kullback-Leibler (KL) divergence
\begin{equation}\label{eq:dKL}
d_{\tKL}(p||\hat{p}):=\int_{\mR^n} p\log \left(\frac{p}{\hat{p}}\right) dx
\end{equation}
represents a well-accepted notion of distance between the two
\cite{KL_informaton,cover2008elements}.
In the case where $p$ and $\hat{p}$ are normal with zero-mean and covariances $T$ and $\hat T$, respectively, their KL divergence becomes
\begin{align*}
 d_{\tKL}(p|| \hat{p})&
=\frac12 \left( \log |\hat{T}|-\log |T|+\trace(T\hat{T}^{-1})-n \right), \nonumber
\end{align*}
while
\begin{align*}
 d_{\tKL}(\hat{p}||p )&=\int_{\mR^n} \hat{p}\log \left(\frac{\hat{p}}{p}\right)dx \nonumber\\
 &=\frac{1}{2}(-\log|\hat{T} |+\log|T|+\trace (\hat{T} T^{-1})-n)\nonumber\\
 &=d_{\tL}(T||\hat T).
 \end{align*}

\subsection{Fisher metric and geodesic distance}

The KL divergence induces a Riemannian structure on the manifold of probability distributions. The quadratic term of $d_{\tKL}(\hat{p}||p+\delta )$ in the perturbation $\delta$ is the Fisher information metric
\begin{align}\label{metric_Fisher}
\g_{p, \text{Fisher}}(\delta)= \int \frac{\delta^2}{p}dx.
\end{align}
This turns out to be natural from one additional perspective. It is the unique Riemannian metric for which the stochastic maps are contractive \cite{Cencov_statistical} --a property that motivates a rich family of metrics in the context of matricial counterparts of probability distributions (see below).

For probability distributions $p(x,\theta)$ parameterized by a vector $\theta$ the corresponding metric is often referred to as Fisher-Rao \cite{Amari_differential} and given by
\[
\g_{p, \text{Fisher-Rao}}(\delta_{\theta})=\delta_{\theta}'E\left[ \left( \frac{\partial \log p}{\partial\theta}\right)\left( \frac{\partial \log p}{\partial\theta}\right)'  \right]\delta_{\theta}.
\]
For zero-mean Gaussian distributions parameterized by corresponding covariance matrices the metric
becomes
\begin{align}\label{metric_Rao}
\g_{T, \text{Rao}}(\Delta)=\left\|T^{-1/2} \Delta T^{-1/2} \right\|_F^2
\end{align}
and is often named after C.R.\ Rao.  We summarize this below. Throughout $\|M\|_F$ denotes the Frobenius norm $\sqrt{{\rm trace}(MM^\prime)}$.

\begin{prop}\label{prop_MetricEqui}
Consider a zero-mean, normal distribution $p$ with covariance $T>0$, and a perturbation $p_\epsilon$ with covariance $T+\epsilon\Delta$. Provided $|| T^{-1/2}\epsilon\Delta T^{-1/2}||_F<1$,
\[
d_{\tKL}(p|| p_\epsilon)= \frac{1}{4}\g_{T, \text{Rao}}(\epsilon\Delta)+O(\epsilon^3).
\]
Moreover, for $\delta=p_\epsilon-p$,
\[\g_{T, \text{Rao}}(\epsilon\Delta)=2 \g_{p, \text{Fisher}}(\delta)+O(\epsilon^4).
\]
\end{prop}
\vspace*{.1in}
The proof is given in Appendix \ref{proof_MetricEqui}.

The Fisher-Rao metric has been studied extensively in recent years \cite{Amari_differential,Bhatia_positive}.
Geodesics and the geodesic distance on the respective Riemannian manifolds can be computed explicitly. In fact, on the space of covariance matrices, the geodesic distance
between two points $T$ and $\hat{T}$ is precisely the log-deviation
 \[
 d_{\tLog}(T, \hat{T})=\|\log (\hat{T}^{-1/2}T \hat{T}^{-1/2}) \|_F.
 \]
Two properties that are worth noting is that the metric is congruence invariant and that the corresponding metric space is complete.

\subsection{Bures metric and Bures/Hellinger distance}
As noted earlier, the Fisher information metric is the unique Riemannian metric for which stochastic maps are contractive. In quantum mechanics, a similar property has been sought for the non-commutative analog of probability vectors, namely, density matrices. These are positive semi-definite and have trace equal to one. In this setting, there are several metrics for which stochastic maps (these are now linear maps between spaces of density matrices, preserving positivity and trace) are contractive. They take the form
\[
\trace(\Delta D_{T}(\Delta))
\]
where $D_T(\Delta)$ can be thought of as a ``non-commutative division'' of the matrix $\Delta$ by the matrix $T$. Thus, if $T,\Delta$ are scalars, the above collapses to $\Delta^2/T$. Particular expressions generating such a ``non-commutative division'' are
\begin{subequations}
\begin{align}
D_{T,1}(\Delta)&:=T^{-1}\Delta,\label{metric1}\\
D_{T,2}(\Delta)&:=\int_{0}^\infty (T+sI)^{-1}\Delta (T+sI)^{-1}ds,\label{metric2}\\
D_{T,3}(\Delta)&:=M, \text{~where~}\frac12 (TM+MT)=\Delta, \label{metric3}
\end{align}
\end{subequations}
see e.g., \cite{Petz1994geometry}.
The metric corresponding to \eqref{metric1} was studied by Petz \cite{Petz1994geometry}, the metric corresponding to \eqref{metric2} is induced by the von Neumann entropy on density matrices and is known as the Kubo-Mori metric \cite{Petz1994geometry}, while
 \eqref{metric3} gives rise to the Bures metric.

The Bures metric can also be written as
 \[
 \g_{T, \text{Bures}}(\Delta):=\min_W \{~\|Y\|_F^2 \mid \Delta=YW'+WY', T=WW' \},
 \]
see \cite{Uhlmann_metric}. Accordingly, the corresponding geodesic distance on the manifold of density matrices is called the Bures length.
Assuming the normalization ${\rm trace}(T)=1$ (i.e., that $T\geq 0$ is a density matrix), $\|W\|_F^2=1$. Thus, we can regard $W$ as an element on a unit sphere. Then, the Bures length is the arc length between corresponding points on the sphere.

The Bures metric has a close connection to the so-called Hellinger distance. A generalization of the standard Hellinger distance to matrices proposed in Ferrante {\em etal.} \cite{Ferrante_hellinger} is
\begin{align}
d_{\tH}(T, \hat{T}):&=\min_{U, V} \left\{\|T^{\frac12}U-\hat{T}^{\frac12}V \|_F\mid UU'=I,VV'=I\right\}\nonumber\\
&=\min_{U} \left\{\|T^{\frac12}U-\hat{T}^{\frac12} \|_F \mid UU'=I\right\},\label{Hellinger_definition}
\end{align}
since, clearly, only one unitary transformation $U$ can attain the same minimal value. This differs from the more standard way to generalize the scalar Hellinger distance to matrices  which is $\trace((T^{1/2}-\hat{T}^{1/2})^2)$. The generalization in (\ref{Hellinger_definition}) is better known in quantum mechanics literature as the Bures distance when the matrices are normalized to have trace $1$.
It is seen that the Bures/Hellinger distance represents a ``straight line'' distance between representatives of two ``points'' on the sphere as measured when imbedded in a linear Euclidean space. The representatives amount to a selection of suitable points on an equivalence class defined via unitary transformations.

Interestingly, as shown in \cite{Uhlmann_metric,Ferrante_hellinger},
\[
d_{\tH}(T, \hat{T})=\left(\trace(T+\hat{T}-2(\hat{T}^\frac12 T \hat{T}^\frac12)^\frac12) \right)^\frac12.
\]
Also, the optimizing unitary matrix $U$ in \eqref{Hellinger_definition} is
\[
U=T^{-\frac12}\hat{T}^{-\frac12}(\hat{T}^\frac12 T \hat{T}^\frac12)^\frac12.
\]
We note that the Hellinger distance applies equally well to positive definite matrices without any need to normalize $T$ and $\hat{T}$, and as such, it
has been used to compare multivariate power spectral densities \cite{Ferrante_hellinger}.

\subsection{Transportation distance}\label{subsect:transportation}

We shift to a seemingly different way of comparing pdf's.
The transportation distance quantifies the cost for transferring one
``mass'' distribution to another accounting for the combined cost of moving every unit of mass from one location to another. Background on transportation problems goes back to the work of G.\ Monge in the 1700's. The recent interest was sparked by the developments in the 1940's by L.\ Kantorovich who is considered the father of the subject\footnote{L.\ Kantorovich received the Nobel prize in Economics in 1975 for his related work on mass transport and resource allocation.}. The importance of transportation distances in probability theory stems from the fact that the respective metrics are weakly continuous.

We consider distributions in ${\mathbb R}^n$ and a quadratic cost. A formulation of the Monge-Kantorovich transportation problem (with a quadratic cost) directly in probabilistic terms is as follows.
Let $X$ and $Y$ be random variables in ${\mathbb R}^n$ having pdf's $p_x$ and $p_y$.
Determine
\begin{align}\label{transportation_define}
d^2_{\rm{W_2}}(p_x,p_y):=\inf_{p}\left\{E (|X-Y|^{2})\mid
\int_x p=p_y,\int_y p=p_x \right\}.
\end{align}
The metric $d_{W_2}$ is known as the Wasserstein metric and, quite surprisingly, also induces a Riemannian structure on probability densities \cite{benamou2000computational,jordan1998variational} -- a rather deep result. Returning to the above optimization, the cost is simply the minimum variance when the marginals of the joint distribution are specified.

We now assume that $T$ and $\hat T$ are the covariances of $X$ and $Y$, respectively, and we let $S=E(XY^\prime)$ denote their correlation. Further, assuming that their joint distribution is Gaussian we obtain
\begin{equation}
\begin{split}
d^2_{\rm{W_2}}(p, \hat{p})= &\min_{S}\left\{ \trace(T+\hat{T}-S-S') \mid \right.\\
  & \left. \hspace*{2.5cm}\left[
          \begin{array}{cc}
         T & S \\
         S' & \hat{T} \\
         \end{array}
         \right]\geq 0
 \right\}.
\end{split}
\end{equation}

A closed form solution is easy to obtain \cite{Olkin_distance,Knott_optimal}:
\begin{align}\label{S0expression}
S_0=\hat{T}^{-\frac12} (\hat{T}^{\frac12} T \hat{T}^{\frac12})^{\frac12}\hat{T}^{\frac12},
\end{align}
and the transportation distance is given alternatively by
\[
d_{\rm{W_2}}(p, \hat{p})=\left(\trace(T+\hat{T}-2(\hat{T}^\frac12 T \hat{T}^\frac12)^\frac12) \right)^\frac12.
\]
Since this is central to our theme, we provide details in Appendix \ref{proof_transp}.
Comparing now with the corresponding expression for the Hellinger distance we readily have the following.
\begin{prop}\label{prop:HellingerTransport}
For $p$ and $\hat{p}$ Gaussian zero mean distributions with covariances $T$ and $\hat T$, respectively,
\[
d_{\tH}(T,\hat{T})=d_{\rm{W_2}}(p,\hat{p}).
\]
\end{prop}

\section{Approximation of structured covariances}\label{sec:approximation}
Returning to the structured covariance approximation problem, we consider the computation of the optimizers for \eqref{approximation}. We do this for every choice of distance discussed in the previous section.

\subsection{Approximation based on KL divergence and likelihood}
If we use $d_{\tKL}$ given in \eqref{eq:dKL} as the distance between $T$ and $\hat{T}$, for the approximation problem we need to solve
\begin{align}\label{problem:KL}
\min_{T\in \cT}\left\{\log |\hat{T}|-\log |T|+\trace(T\hat{T}^{-1})-n\right\}.
\end{align}
This is convex in $T$, provided $\hat{T}>0$, and hence numerically feasible. However, the problem is vacuous when $\hat{T}$ is singular. This is unsatisfactory since the case when $\hat T$ is singular is important and quite common.
Alternatively, if we
use the likelihood divergence $d_{\tL}(T||\hat{T} )$ as distance measure, the optimization problem
\begin{align}\label{problem:maximum}
\min_{T\in \cT}\left\{-\log|\hat{T}|+\log|T|+\trace (\hat{T} T^{-1})-n\right\}
\end{align}
is well defined for singular $\hat T$ as well.

A necessary condition for a local minimum of \eqref{problem:maximum} given in \cite{BLW_estimation} is:
\begin{align}\label{necessary}
\trace\left((T^{-1}\hat{T} T^{-1}-T^{-1})Q\right)=0,
\end{align}
for all Toeplitz $Q$ and pointed out in \cite{BLW_estimation} that, provided $\hat{T} $ is not singular, there is at least one local minimum of \eqref{problem:maximum} which is positive definite. Based on this, Burg {\em etal} \cite{BLW_estimation} give a numerical method to solve \eqref{problem:maximum}. The method is computationally demanding and numerically sensitive, especially when $\hat T$ is singular.

\subsection{Approximation based on log-deviation}

The optimization problem
\begin{align}\label{problem:Log}
\min_{T\in \cT}\left\{ \|\log (\hat{T}^{-1/2}T \hat{T}^{-1/2}) \|_F\right\}.
\end{align}
is not convex in $T$. Linearization of the objective function about $\hat{T}$ may be used instead, since this leads to
\begin{align}\label{problem:LogLinear}
\min_{T\in \cT}\left\{ \|\hat{T}^{-1/2}T \hat{T}^{-1/2}-I \|_F\right\}
\end{align}
which is a convex problem.

\subsection{Based on Hellinger and transportation distance}

Using $d_{\tH}(T,\hat{T})$ the relevant optimization problem \eqref{approximation} becomes
\begin{align}\label{problem:Hellinger}
\min_{T\in\cT}\left\{\trace(T+\hat{T}-2(\hat{T}^\frac12 T \hat{T}^\frac12)^\frac12)\right\}.
\end{align}
At the outset, this appears difficult. However, from Proposition \ref{prop:HellingerTransport} we know that $d_{\tH}(T,\hat{T})=d_{\rm{W_2}}(p,\hat{p})$. Hence, we may evaluate \eqref{problem:Hellinger} via solving
\begin{align}\label{problem:transp}
\min_{T\in \cT,\; S}\left\{\trace(T+\hat{T}-S-S') \mid
 \left[                                                                 \begin{array}{cc}
                                                                                           T & S \\
                                                                                           S' & \hat{T} \\
                                                                                         \end{array}
                                                                                    \right]\geq 0\right\}   .
\end{align}
This is now a semi-definite program and can be solved quite efficiently \cite{Boyd}.

The above expression for the transportation distance can be given an {\em alternative
interpretation} as follows. We postulate the statistical model
\newcommand{\bv}{{\mathbf v}}
\[
\hat{\bx}=\bx+\bv
\]
where $\bv$ represents noise, and $\hat T$ and $T$ are the covariances of
$\hat \bx$ and $\bx$, respectively. The covariance of $\bx$ is known to be in the admissible set $\cT$
while that of $\hat \bx$ may not, due to noise. Thus, in the absence of additional priors, it is reasonable to seek an ``explanation'' of the estimated covariance $\hat T$ by assuming the least possible amount of noise.
Allowing for possible coupling between $\bx$ and $\bv$ brings us to minimize
\[
E\{\bv^\prime\bv\}=\trace(T+\hat{T}-S-S')
\]
subject to positive semi-definiteness of the covariance of $[\bx^\prime,\;\hat\bx^\prime]^\prime$. This is precisely \eqref{problem:transp}.

Analogous rationale, albeit with different assumptions has been used to justify different methods.
For instance, assuming that $\hat\bx=\bx+\bv$ where $\bx$ and $\bv$ are independent leads to
\begin{align*}
\min_{T\in\cT}\left\{ \trace(\hat{T}-T)\mid \hat{T}-T\geq0\right\}
\end{align*}
which is a method proposed in \cite{Stoica_MA}. Then, also,
assuming a ``symmetric'' noise contribution as in
\[
\hat \bx+\hat \bv=\bx+\bv,
\]
where the noise vectors $\hat{\bv}$ and $\bv$ are independent of $\bx$ and $\hat \bx$, leads
to
\begin{align*}
\min_{T\in\cT, Q,\hat{Q}}\left\{ \trace(\hat{Q}+Q)\mid \hat{T}+\hat{Q}=T+Q,~Q,\hat{Q}\geq0\right\},
\end{align*}
where $\hat{Q}$ and $Q$ designate
covariances of $\hat{\bv}$ and $\bv$,  respectively. The minimum in this case is the
nuclear norm of $\hat{T}-T$ and studied as a possibility in \cite{Georgiou_L1Distances}.

\section{Examples}\label{sec:examples}
We now compare how well two of the methods outlined earlier perform in identifying a single spectral line in white noise and compare those with the standard Burg's method. We choose parameters as in the example in Burg {\em etal.}  \cite{BLW_estimation}. For constructing power spectra corresponding to a finite set of covariance samples and we use autoregressive models in order to be consistent with \cite{BLW_estimation}. By using the same type of power spectra we isolate and compare the effect of correcting for the ``non-Toeplitz-ness'' via each of these two methods and by Burg's method.

The data consists of a sinusoid (leading to a single spectral line) with three different phase values and the same random vector for the noise. We assume a single observation vector of size $11$, hence both $\bx$ and $\bv$ are vectors, and thus, the estimated covariance $\hat T$ is $(11\times 11)$, singular, and of rank equal to $1$.
Thus, the data is the same additive mixture of sinusoid and noise as in \cite{BLW_estimation}:
\[
x(t)=\cos(\frac{\pi}{4}t+\psi)+v(t), t=0,1,\ldots,10.
\]
The initial phase $\psi$ is chosen for three different values $\frac{\pi}{4},\frac{\pi}{2}$ and $\frac{3\pi}{4}$.
The noise vector $v$ is fixed as
\begin{align*}
v=[&0.000562,~ -0.019127,~ 0.007377,~ -0.000149,  -0.007479,~-0.013960,\\
&0.003510, ~~~~0.012380, ~~0.006979, ~~~~0.003092,~~ 0.010053]',
\end{align*}
generated from a zero-mean Gaussian distribution with variance $0.0001$.

\begin{figure}[htb]\begin{center}
\includegraphics[totalheight=9cm]{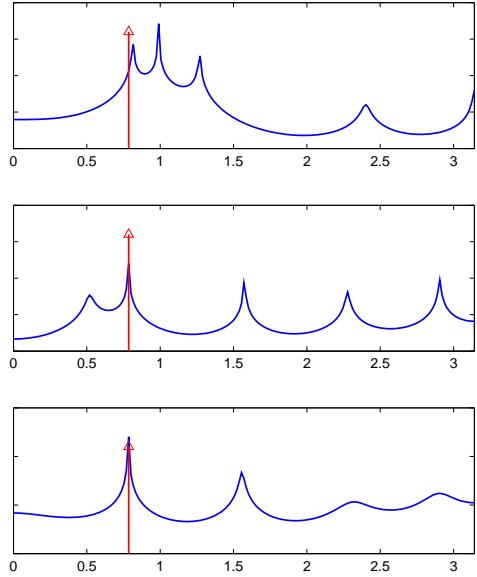}
\caption{Estimated maximum entropy spectrum for $\psi=\frac{\pi}{4}$: i) Burg's method, ii) Maximum likelihood method, iii) Minimum transportation method.}\label{fig:example1pi4}\end{center}
\end{figure}

The first plot in Figure \ref{fig:example1pi4} shows the power spectral density (PSD) using Burg's method for estimating the partial correlation coefficients (as in \cite{BLW_estimation}), while the second and third plots are based on covariances approximated using the likelihood-based method and transportation-based methods, respectively. The data corresponds to $\psi=\frac{\pi}{4}$ and the resolution of the plots is $\frac{\pi}{200}$. The (red) arrow in the plots indicates the frequency of the sinusoidal component. Burg's method splits the spectral line into three. The spectral line closest to the true (red arrow) is also significantly off. On the other hand, both, the likelihood-based and the transportation-based methods detect the spectral line at the correct frequency (with relatively insignificant error).

\begin{figure}[htb]\begin{center}
\includegraphics[totalheight=9cm]{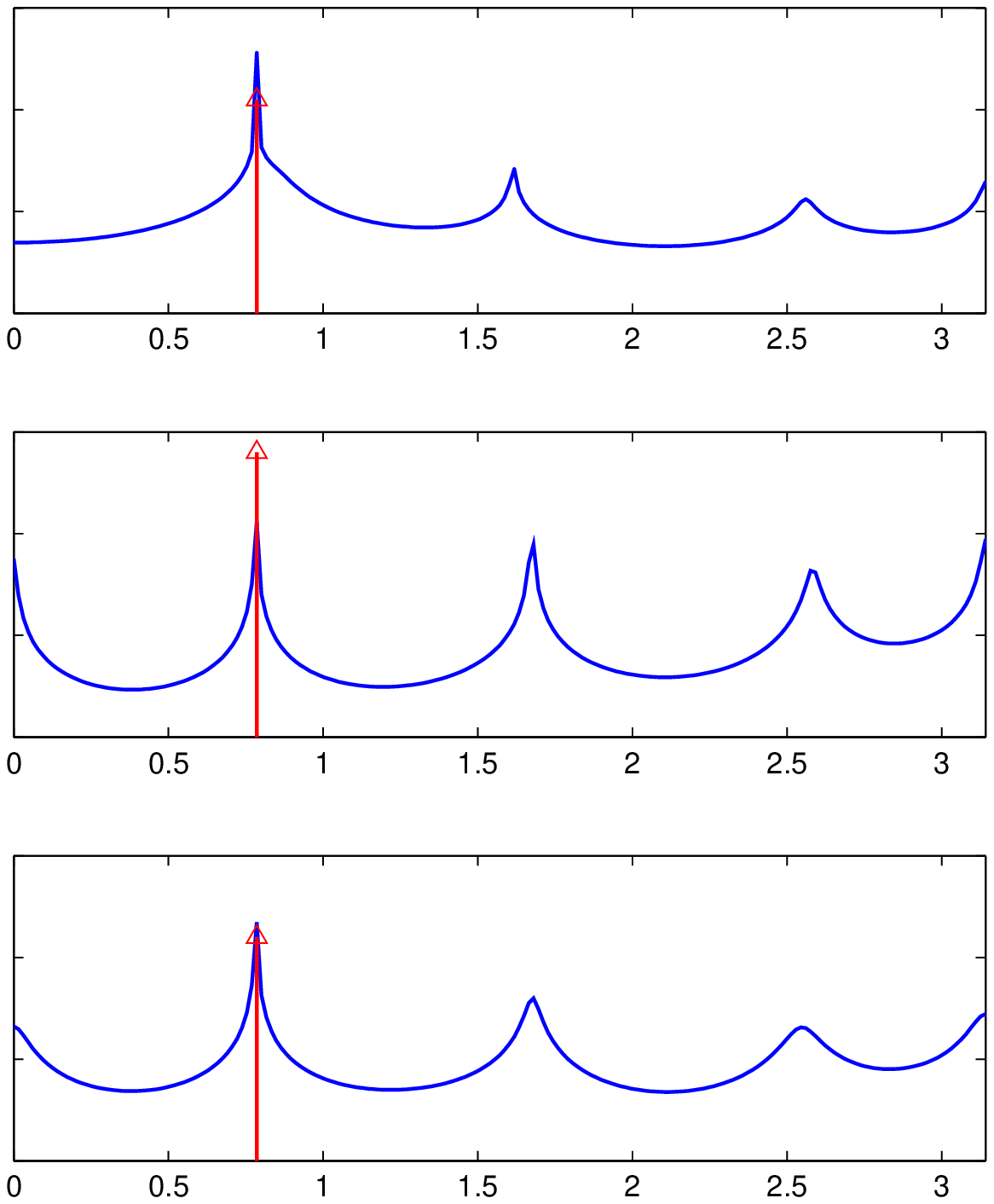}
\caption{Estimated maximum entropy spectrum for $\psi=\frac{\pi}{2}$: i) Burg's method, ii) Maximum likelihood method, iii) Minimum transportation method.}\label{fig:example1pi2}\end{center}
\end{figure}

\begin{figure}[htb]\begin{center}
\includegraphics[totalheight=9cm]{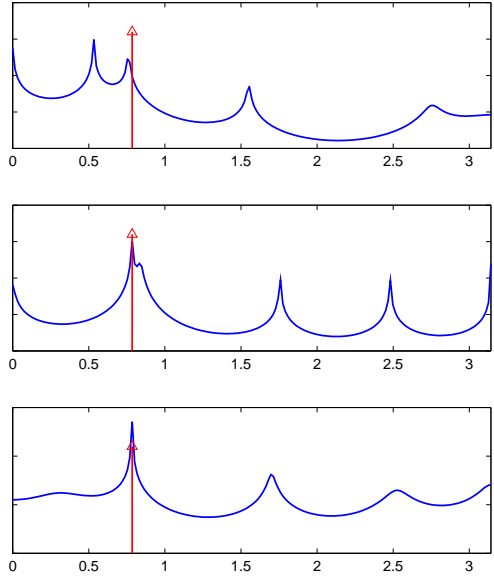}
\caption{Estimated maximum entropy spectrum for $\psi=\frac{3\pi}{4}$: i) Burg's method, ii) Maximum likelihood method, iii) Minimum transportation method.}\label{fig:example13pi4}\end{center}
\end{figure}

Figure \ref{fig:example1pi2} shows the same situation but for $\psi=\frac{\pi}{2}$. All three methods detect the spectral line perfectly, to within the stated resolution.  Figure \ref{fig:example13pi4} corresponds to the case where $\psi=\frac{3\pi}{4}$. Burg's method consistently splits the true spectral line into two nearby ones. The likelihood-based method gives a small peak near the true spectral line, although the dominant line is located at the true frequency to within the stated resolution. On the other hand, the transportation-based method gives a result which is consistent with the previous two situations. For the purpose of detecting line spectra, the transportation-based method appears to be the most robust.

A potential drawback of the transportation-based method is that it gives a biased estimate for the energy in the sinusoidal component. This is typically smaller than the true value in the example.


\section{Recap} \label{sec:conclusions}

Most modern spectral analysis methods rely on estimated covariance statistics.
Yet, they are sensitive to those statistics abiding by the requisite linear structure, e.g., Toeplitz.
In this paper we discussed and compared two of the most promising methods
for approximating a sample covariance with one of the required structure. Contributions in the paper
include drawing the connection between approximation in the Hellinger distance and approximation in the sense of optimal mass transport. The latter can be cast as a semidefinite program which is easy to solve and impervious to possible singularity or near-singularity of the sample covariance.

The issue with the sample covariance $\hat{T}$ being singular is often neglected in estimation problems. Yet, it is ubiquitous when only few short observation records are available ---a situation which is common in the analysis of non-stationary processes.
Furthermore, the uniqueness and other properties of a maximum likelihood estimate, when $\hat T$ is singular, are not well understood \cite{BLW_estimation}.

As a final remark we note that interest in other linear structures for covariance matrices, besides that of Toeplitz, arises when the vectorial process is the state vector of a linear system. In such a case,
$T$ satisfies linear constraints that involve the system dynamics \cite{Georgiou_Statecovariance}. All earlier discussion and methods can be repeated verbatim for the problem of approximating sample state-covariances.

\section{Appendix A: proof for the Proposition \ref{prop_MetricEqui}}\label{proof_MetricEqui}

The KL divergence between a zero-mean normal distribution $p$ with covariance $T>0$ and a perturbation $p_\epsilon$ with covariance $T+\epsilon\Delta$ is
\begin{align*}
d_{\tKL}(p||p_\epsilon)=\frac{1}{2}\left(\log{\det(T+\epsilon\Delta)}-\log{\det(T)}+\trace\left((T+\epsilon\Delta)^{-1}T\right)-n  \right).
\end{align*}
\newcommand{\DeltaT}{{\Delta_{T}}}
Define $\DeltaT =T^{-1/2}\Delta T^{-1/2}$, then
\begin{align}\label{CM_5}
d_{\tKL}(p||p_\epsilon)&=\frac{1}{2}\left(\log{\det\left(T^{1/2}(I+\epsilon \DeltaT )T^{1/2}\right)}-\log{\det(T)}+\trace\left(T^{-1/2}(I+\epsilon \DeltaT )^{-1}T^{-1/2}T\right)-n\right)\nonumber\\
&=\frac{1}{2}\left(\log{\det(I+\epsilon \DeltaT )}+\trace(I+\epsilon \DeltaT )^{-1}-n  \right).
\end{align}
We expand $(I+\epsilon \DeltaT )^{-1}$ into the Taylor series
\begin{align}\label{CM_6}
(I+\epsilon \DeltaT )^{-1}=I-\epsilon\DeltaT +\epsilon^2 \DeltaT ^2-\epsilon^3 \DeltaT ^3+\cdots.
\end{align}
Let $\lambda_i, ~i=1, ~\cdots, ~n$, represent eigenvalues of $\DeltaT $, then
\begin{align}\label{CM_7}
\log{\det(I+\epsilon \DeltaT )}&=\sum_{i=1}^n \log(1+\epsilon \lambda_i)\nonumber\\
&=\sum_{i=1}^n (\epsilon \lambda_i-\frac12 \epsilon^2 \lambda_i^2+\frac{1}{3}\epsilon^3 \lambda_i^3+\cdots)\nonumber\\
&=\epsilon \trace(\DeltaT )-\frac12 \epsilon^2 \trace(\DeltaT ^2)+\frac{1}{3}\epsilon^3 \trace(\DeltaT ^3)+\cdots.
\end{align}
We substitute (\ref{CM_6}) and (\ref{CM_7}) into (\ref{CM_5}) to obtain
\[
d_{\tKL}(p||p_\epsilon)=\frac{1}{4}\epsilon^2 \trace(\DeltaT ^2)+O(\epsilon^3).
\]
By a similar computation, one can easily see that $d_{\tKL}(p_\epsilon||p)$ gives rise to the same metric,
though the coefficients of higher order terms on $\epsilon$ are different from those corresponding to $d_{\tKL}(p||p_\epsilon)$.

To draw a connection with the Fisher metric, we substitute $\delta=p_\epsilon-p$ into the Fisher metric:
\begin{align*}
\g_{p, \text{Fisher}}(\delta)=\left( \int_{\mR^n} \frac{\det(T)^{1/2}}{(2\pi)^{n/2}\det(T+\epsilon \Delta)}e^{-\frac12 y'(2(T+ \epsilon \Delta)^{-1}-T^{-1})y}dy-1\right).
\end{align*}
Since $||\epsilon \DeltaT  ||_F<1$, $\epsilon^2\DeltaT ^2< I$ and hence
\[
-I<\epsilon \DeltaT < I.
\]
Multiplying by $T^{1/2}$ from left and right  on all sides of the above inequality, we obtain
\[
-T<\epsilon {\Delta}<T,
\]
or equivalently
\[
0<\frac12 T+\frac12 \epsilon {\Delta}<T.
\]
It follows that
\[
-(\frac12 T+\frac12 \epsilon {\Delta} )^{-1}< -T^{-1},
\]
or equivalently,
\[
2(T+\epsilon {\Delta} )^{-1}-T^{-1}>0.
\]
Consequently
\[
\frac{1}{(2\pi)^{n/2}\det\left(2(T+\epsilon {\Delta} )^{-1}-T^{-1}\right)^{-1/2}}e^{-\frac12 y'(2(T+ \epsilon \Delta)^{-1}-T^{-1})y}
\]
is a Gaussian distribution with mean $0$ and covariance $(2(T+ \epsilon \Delta)^{-1}-T^{-1})^{-1}$. Since the integral of a Gaussian distribution is $1$, we obtain that
\[
\g_{p, \text{Fisher}}(\delta)=\left( \frac{\det(T)^{1/2}}{\det\left(2(T+\epsilon \Delta)^{-1}-T^{-1}\right)^{1/2} \det(T+\epsilon \Delta)} -1 \right).
\]
But
\begin{align*}
(T+\epsilon \Delta)^{-1}&=T^{-1/2}(I+\epsilon \DeltaT )^{-1} T^{-1/2},
\end{align*}
and
\[
2(T+\epsilon \Delta)^{-1}-T^{-1}=T^{-1/2}\left(2(I+\epsilon \DeltaT )^{-1}-I\right)T^{-1/2}.
\]
Consequently,
\begin{align*}
\det\left(2(T+\epsilon \Delta)^{-1}-T^{-1}\right)^{1/2}\det(T+\epsilon \Delta)&=\det\left( (T+\epsilon \Delta)\left(2(T+\epsilon \Delta)^{-1}-T^{-1}\right)(T+\epsilon \Delta)\right)^{1/2}\\
&=\det\left(T^{1/2}(I+\epsilon \DeltaT )\left(2(I+\epsilon \DeltaT )^{-1}-I\right)(I+\epsilon \DeltaT )T^{1/2}\right)^{1/2}\\
&=\det(T)^{1/2}\det(I-\epsilon^2 \DeltaT ^2)^{1/2},
\end{align*}
and
\begin{align*}
\g_{p, \text{Fisher},m}(\delta)=\left(\det(I-\epsilon^2\DeltaT ^2)^{-1/2}-1\right)=\left(\det(I+\epsilon^2\DeltaT ^2+\epsilon^4\DeltaT ^4+\cdots)^{1/2}-1\right).
\end{align*}
Once again considering the eigenvalues of $\DeltaT$ we get
\begin{align*}
\det(I+\epsilon^2\DeltaT ^2+\epsilon^4\DeltaT ^4+\cdots)^{1/2}&=\left(\prod_{k=1}^n (\sum\limits_{i=0}^{\infty}(\epsilon\lambda_k)^{2i})\right)^{1/2}\\
&=\left(1+\sum\limits_{k=1}^n \epsilon^2\lambda_k^2+\sum\limits_{k\leq l}\epsilon^4 \lambda_k^2\lambda_l^2+\cdots \right)^{1/2}\\
&=1+\frac12\epsilon^2\|\DeltaT \|_F^2+O(\epsilon^4),
\end{align*}
where in the last equality we have used the fact that $\sum\limits_{k=1}^n \lambda_k^2=\trace(\DeltaT ^2)=\|\DeltaT \|_F^2$. Therefore,
\[
\g_{p, \text{Fisher}}(\delta)=\frac12 \g_{T, \text{Rao} }(\Delta)+O(\epsilon^4).
\]

\section{Appendix B}\label{proof_transp}
We now show that given two $n\times n$ matrices $T> 0$ and $\hat{T}> 0$,
\[
{\rm arg}\min_{S}~\left\{\trace (T+\hat{T}-S-S') \mid \left[
                           \begin{array}{cc}
                             T & S \\
                             S' & \hat{T} \\
                           \end{array}
                         \right]\geq0
 \right\}
\]\\
has indeed the explicit closed-form expression
\begin{align}
S_0=\hat{T}^{-\frac12}(\hat{T}^\frac12 T \hat{T}^\frac12)^\frac12 \hat{T}^\frac12, \label{lemma1:eq1}
\end{align}

Consider the Shur complement
\[
P:=T-S\hat{T}^{-1}S'
\]
which is clearly nonnegative definite. Then, $S\hat{T}^{-\frac12}=(T-P)^\frac12U$, where $UU'=I$, and
\begin{align}\label{eq:C}
S=(T-P)^\frac12 U \hat{T}^\frac12.
\end{align}
Moreover, \begin{align}\label{traceC}
\trace(S)=\trace((T-P)^\frac12 U \hat{T}^\frac12)= \trace (\hat{T}^\frac12(T-P)^\frac12 U).
\end{align}
Since $T$ and $\hat{T}$ are given, minimizing $\trace (T+\hat{T}-S-S')$ is the same as maximizing
$\trace(S)$.
Let $U_S \Lambda_S V_S'$ be the singular value decomposition of
$
\hat{T}^\frac12(T-P)^\frac12$,  and
\[
U_0:=\arg\max_U~ \{\trace (\hat{T}^\frac12(T-P)^\frac12 U) \mid UU'=I\}.
\]
Then, $U_0$ must satisfy $V_S'U_0=U_S'$ and
\begin{align}\label{symmetric}
 \hat{T}^\frac12(T-P)^\frac12 U_0=(\hat{T}^\frac12(T-P)\hat{T}^\frac12)^\frac12.
\end{align}
From \eqref{traceC} we have
$
\trace (S)=\trace( (\hat{T}^\frac12(T-P)\hat{T}^\frac12)^\frac12).
$
Since $P\geq0$, the $\trace (S)$ is maximal when $P=0$.
Moreover, if $P=0$,
 \[
 \rank\left(
 \left[
                           \begin{array}{cc}
                             T & S \\
                             S' & \hat{T} \\
                           \end{array}
                         \right]\right)\leq \rank(T),
 \]
 and $
 \hat{T}=S_0'T^{-1}S_0.
$
Thus, setting $P=0$ into \eqref{symmetric}, we have
\[
U_0=T^{-\frac12}\hat{T}^{-\frac12}(\hat{T}^\frac12 T \hat{T}^\frac12)^\frac12,
\]
and consequently
$
S_0=\hat{T}^{-\frac12}(\hat{T}^\frac12 T \hat{T}^\frac12)^\frac12\hat{T}^\frac12.
$

\bibliographystyle{IEEEtran}
\bibliography{IEEEabrv,CovEst}

\begin{thebibliography}{10}
\providecommand{\url}[1]{#1}
\csname url@samestyle\endcsname
\providecommand{\newblock}{\relax}
\providecommand{\bibinfo}[2]{#2}
\providecommand{\BIBentrySTDinterwordspacing}{\spaceskip=0pt\relax}
\providecommand{\BIBentryALTinterwordstretchfactor}{4}
\providecommand{\BIBentryALTinterwordspacing}{\spaceskip=\fontdimen2\font plus
\BIBentryALTinterwordstretchfactor\fontdimen3\font minus
  \fontdimen4\font\relax}
\providecommand{\BIBforeignlanguage}[2]{{%
\expandafter\ifx\csname l@#1\endcsname\relax
\typeout{** WARNING: IEEEtran.bst: No hyphenation pattern has been}%
\typeout{** loaded for the language `#1'. Using the pattern for}%
\typeout{** the default language instead.}%
\else
\language=\csname l@#1\endcsname
\fi
#2}}
\providecommand{\BIBdecl}{\relax}
\BIBdecl

\bibitem{grenander2001toeplitz}
U.~Grenander and G.~Szeg{\"o}, \emph{Toeplitz Forms and their
  Applications}.\hskip 1em plus 0.5em minus 0.4em\relax Chelsea Pub Co, 2001.

\bibitem{Burg_thesis}
J.~Burg, ``Maximum entropy spectral analysis,'' Ph.D. dissertation, Stanford
  University, Stanford, CA, 1975.

\bibitem{haykin1979nonlinear}
S.~Haykin, \emph{Nonlinear Methods of Spectral Analysis}.\hskip 1em plus 0.5em
  minus 0.4em\relax Springer-Verlag, 1979.

\bibitem{BLW_estimation}
J.~Burg, D.~Luenberger, and D.~Wenger, ``Estimation of structured covariance
  matrices,'' \emph{Proceedings of the IEEE}, vol.~70, no.~9, pp. 963--974,
  1982.

\bibitem{KL_informaton}
S.~Kullback and R.~A. Leibler, ``On information and sufficiency,'' \emph{The
  Annals of Mathematical Statistics}, vol.~22, no.~1, pp. 79--86, 1951.

\bibitem{cover2008elements}
T.~Cover and J.~Thomas, \emph{Elements of Information Theory}.\hskip 1em plus
  0.5em minus 0.4em\relax Wiley-Interscience, 2008.

\bibitem{Cencov_statistical}
N.~Cencov, \emph{Statistical Decision Rules and Optimal Inference}.\hskip 1em
  plus 0.5em minus 0.4em\relax Amer. Math. Soc., 1982, no.~53.

\bibitem{Amari_differential}
S.-I. Amari, ``Differential-geometrical methods in statistics,'' \emph{Lecture
  Notes in Statistics}, vol.~28, 1985.

\bibitem{Bhatia_positive}
R.~Bhatia, \emph{{Positive Definite Matrices}}.\hskip 1em plus 0.5em minus
  0.4em\relax Princeton University Press, 2007.

\bibitem{Petz1994geometry}
D.~Petz, ``Geometry of canonical correlation on the state space of a quantum
  system,'' \emph{Journal of Mathematical Physics}, vol.~35, pp. 780--795,
  1994.

\bibitem{Uhlmann_metric}
A.~Uhlmann, ``{The metric of Bures and the geometric phase},'' \emph{Quantum
  Groups and Related Topics}, pp. 267--264, 1992.

\bibitem{Ferrante_hellinger}
A.~Ferrante, M.~Pavon, and F.~Ramponi, ``Hellinger versus {K}ullback--{L}eibler
  multivariable spectrum approximation,'' \emph{Automatic Control, IEEE
  Transactions on}, vol.~53, no.~4, pp. 954--967, 2008.

\bibitem{benamou2000computational}
J.~Benamou and Y.~Brenier, ``A computational fluid mechanics solution to the
  {M}onge-{K}antorovich mass transfer problem,'' \emph{Numerische Mathematik},
  vol.~84, no.~3, pp. 375--393, 2000.

\bibitem{jordan1998variational}
R.~Jordan, D.~Kinderlehrer, and F.~Otto, ``The variational formulation of the
  {F}okker-{P}lanck equation,'' \emph{SIAM journal on mathematical analysis},
  vol.~29, no.~1, pp. 1--17, 1998.

\bibitem{Olkin_distance}
I.~Olkin and F.~Pukelsheim, ``The distance between two random vectors with
  given dispersion matrices,'' \emph{Linear Algebra and its Applications},
  vol.~48, pp. 257--263, 1982.

\bibitem{Knott_optimal}
M.~Knott and C.~S. Smith, ``On the optimal mapping of distributions,''
  \emph{Journal of Optimization Theory and Applications}, vol.~43, no.~1, pp.
  39--49, 1984.

\bibitem{Boyd}
S.~Boyd and L.~Vandenberghe, \emph{Convex Optimization}.\hskip 1em plus 0.5em
  minus 0.4em\relax Cambridge University Press, 2004.

\bibitem{Stoica_MA}
P.~Stoica, L.~Du, J.~Li, and T.~Georgiou, ``A new method for moving-average
  parameter estimation,'' in \emph{Conference Record of the Forty Fourth
  Asilomar Conference on Signals, Systems and Computers}, 2010, pp. 1817--1820.

\bibitem{Georgiou_L1Distances}
T.~Georgiou, ``Distances between time-series and their autocorrelation
  statistics,'' \emph{Modeling, Estimation and Control}, pp. 113--122, 2007.

\bibitem{Georgiou_Statecovariance}
------, ``{The structure of state covariances and its relation to the power
  spectrum of the input},'' \emph{Automatic Control, IEEE Transactions on},
  vol. 47(7), pp. 1056--1066, 2002.

\end{thebibliography}
\end{document}